\documentclass[a4paper,UKenglish,cleveref,autoref,thm-restate]{lipics-v2021}
\hideLIPIcs 
\nolinenumbers 
\usepackage{mathtools, todonotes}
\usepackage{tikz}
\usetikzlibrary{calc, patterns, decorations.pathreplacing}
\usetikzlibrary{shapes, positioning, arrows.meta}
\newtheorem{problem}[theorem]{Problem}
\DeclareMathOperator{\tw}{tw}
\DeclareMathOperator{\dtw}{dtw}
\DeclareMathOperator{\tpw}{tpw}

\newcommand{\eps}{\varepsilon}
 
\newcommand{\mc}[1]{\mathcal{#1}}

\newcommand{\bd}[1]{\partial #1}

\newcommand{\s}[1]{\left(#1\right)}

\newcommand{\set}[1]{\{#1\}}
\newtheorem{obs}[theorem]{Observation}
\newtheorem{cor}[theorem]{Corollary}

\title{Trade-off between spread and width for tree decompositions} %TODO Please add

\author{Hans L. Bodlaender}{Department of Information and Computing Sciences, Utrecht University}{h.l.bodlaender@uu.nl}{ https://orcid.org/
0000-0002-9297-3330}{}

\author{Carla Groenland}{Faculty of Electrical Engineering, Mathematics and Computer Science, Technical University Delft}{c.e.groenland@tudelft.nl}{https://orcid.org/
0000-0002-9878-8750}{Research supported by the Dutch Research Council (NWO, VI.Veni.232.073).}

\author{Sergey Norin}{Department of Mathematics and Statistics, McGill University}{sergey.norin@mcgill.ca}{}{Research supported by an NSERC Discovery grant.}

\author{Neil Rahman}{Department of Mathematics and Statistics, McGill University}{neil.rahman@mcgill.ca}{}{Research supported by an NSERC Discovery grant.}

\authorrunning{H. L. Bodlaender, C. Groenland, S. Norin, N. Rahman} 

\Copyright{Hans L. Bodlaender, Carla Groenland, Sergey Norin, Neil Rahman} 



\category{}

\relatedversion{} 

\acknowledgements{We would like to thank Hugo Jacob and David Wood for useful comments and suggestions.}

\begin{document}

\maketitle

\begin{abstract}
The \emph{spread} of a vertex $v$ in a tree decomposition is the number of bags that contain $v$.
We study the trade-off between spread and width in tree decompositions, answering every open question from Wood [arXiv:2509.01140]. 
\begin{itemize}
    \item Wood asked for the infimum of $c > 0$ such that there exists $c'$ such that each graph $G$ has a tree decomposition of width $c\tw(G)$ in which each vertex $v$ has spread at most $c'(d(v)+1)$. We show that the answer is $3$.
    \item We prove a conjecture of Wood, stating that every tree-decomposition of the $(n \times n)$-grid with width $n$ has a vertex with spread $\Omega(n)$. 
    \item  Finally, we answer the last question of Wood by showing that near-optimal average spread can be achieved simultaneously with width $O(\tw(G))$.
\end{itemize}
\end{abstract}

\section{Introduction}
Tree decompositions have been widely studied in both structural and algorithmic contexts, for example due to their importance to graph minor theory~\cite{ROBERTSON198449} and the fact that many problems can be solved in linear time on graphs of bounded treewidth~\cite{bodlaenderlineartimealg,COURCELLE199012}.  

The main criterion of `optimality' for a tree decomposition is usually its width, or, equivalently, the size of the largest bag in the tree decomposition. 
Recently, Wood~\cite{wood2025treedecompositionssmallwidth} proved several results concerning tree decompositions which have additional nice properties, while not sacrificing too much of the width.
Our work continues on this, answering several of the problems posed by Wood.

Besides the width, our first parameter of interest is the spread. The \textit{spread} of a vertex $v$ in a
tree-decomposition $(T,(B_x)_{x\in V(T)})$ is the number of bags $v$ is in: $|\{x \in V (T):v \in B_x\}|$. Since all neighbours of a vertex $v$ must appear together with $v$ in some bag, when the width is held constant, the spread of $v$ must necessarily grow with its degree.

Wood~\cite{wood2025treedecompositionssmallwidth} posed the following problem.
\begin{problem}
\label{prob:spread_constant_min_width}
    What is the infimum of $c$ such that for some $c'$, every graph $G$ with  treewidth $k$ has a tree-decomposition with width at most $(c + o(1))k$, in which each vertex $v \in V (G)$ has spread at most $c'(d(v) + 1)$?
\end{problem} 
Wood~\cite{wood2025treedecompositionssmallwidth} already proved that the answer is at most 14 (with $c'=1$).
\begin{theorem}[Theorem 2 in~\cite{wood2025treedecompositionssmallwidth}]
    Every graph $G$ with treewidth $k$ has a tree-decomposition with width at most
$14k + 13$, such that each vertex $v \in V (G)$ has spread at most $d(v) + 1$.
\end{theorem}
We show that the value of $c$ of Problem 1 is 3.
\begin{theorem}
    \label{thm:spread_constant_min_width_upper_bound}
    For each $c>3$, there is a constant $c'>0$ such that every graph $G$ with  treewidth $k$ has a tree-decomposition with width at most $ck$, in which each vertex $v \in V (G)$ has spread at most $c'(d(v) + 1)$.
\end{theorem}
\begin{theorem}
    \label{thm:spread_constant_min_width_lower_bound}
    For each $c<3$ and for each $c'>0$, there is a graph $G$ of treewidth $k$ such that in every tree-decomposition with width at most $ck$, some vertex $v \in V (G)$ has spread strictly more than $c'(d(v) + 1)$.
\end{theorem}

Wood also proposed a lower bound example for Problem~\ref{prob:spread_constant_min_width}, conjecturing there are no constants $\varepsilon,c>0$ such that $(n\times n)$-grid has a tree-decomposition of width at most $(2-\varepsilon)n$ and spread at most $c$.
This conjecture is false and in fact for any $m\geq n\geq 1$ and $\varepsilon>0$, the $(n\times m)$-grid admits a tree-decomposition of width at most $(1+\varepsilon)n$ in which each vertex has spread at most $O(1/\varepsilon)$ (see Lemma~\ref{lem:good_spread_width_for_grids}). On the other hand, Wood conjectured that every tree-decomposition of the $(n \times n)$-grid with width $n$, has some vertex with spread $\Omega(n)$. We prove this conjecture with the following stronger result.
\begin{theorem}\label{thm:grid_linear_spread}
There is a constant $\delta>0$ such that for all integers $n, a \geq 1$, every tree-decomposition of the $(n \times n)$-grid of width $n+a-1$ has a vertex with spread at least $\delta n/a$.
\end{theorem}

When storing a tree-decomposition $(T,(B_x)_{x\in V(T)})$ of a graph $G$, the memory required depends on $\sum_{x\in V(T)}|B_x|$. Our second parameter of interest is the average spread $\frac1{|V(G)|}\sum_{x\in V(T)}|B_x|$.
Regarding this, Wood~\cite{wood2025treedecompositionssmallwidth} raised the following problem.
\begin{problem}
\label{prob:average_spread}
What is the infimum of the $c' \in \mathbb{R}$ such that for some $c$ every graph $G$ has a
tree-decomposition of width at most $c\tw(G)$ and average spread at most $c'$?
\end{problem}
Wood already proved the upper bound 3.
We show that the answer to this problem is 1.
\begin{theorem}
    \label{thm:average_spread_inf_is_1}
Let $c'>1$. Then there is a constant $c>0$ such that every graph $G$ has a
tree-decomposition of width at most $c\tw(G)$ and average spread at most $c'$.

\end{theorem}

\paragraph*{Related work: domino treewidth, tree-partition-width and persistence}
Bienstock~\cite{Bienstock90} asked whether it is possible for each graph $G$ to find a tree decomposition $(T,(B_t)_{t\in V(T)})$ in which $T$ is a binary tree and the width as well as the maximum diameter of $T_v$ (the subtree of $T$ induced on the bags containing $v$, for $v\in V(G)$), are bounded by a constant depending only on $\tw(G)$ and $\Delta(G)$. 
Ding and Oporowski~\cite{DingO95} answered this question positively in a stronger fashion, providing the first bounds for spread that we are aware of.
\begin{theorem}[Theorem 1.3 in~\cite{DingO95}]
Every graph $G$ has a tree decomposition of width at most $$2^{\tw(G)+1}(\tw(G)+1)-1$$ such that each vertex $v\in V(G)$ has spread at most $$1+2d(v)3^{2^{\tw(G)}}.$$
\end{theorem}
Ding and Oporowski~\cite{DingO95} also studied a special case of tree decompositions with bounded spread, namely those where each vertex has spread at most~2.
They showed that each graph admits such a tree decomposition of width at most $O(\text{tw}(G)^2\Delta(G)^3)$, where $\tw(G)$ and $\Delta(G)$ denote the treewidth and maximum degree of the graph $G$ respectively.

Bodlaender and Engelfriet~\cite{BodlaenderE97} first used the term \textit{domino treewidth} for the optimal width of a tree decomposition in which each vertex has spread 2. We write $\dtw(G)$ for the domino treewidth of $G$. Bodlaender~\cite{Bodlaender99} improved the upper bound on domino treewidth to $\dtw(G)=O(\tw(G)\Delta(G)^2)$ and gave examples of graphs $G$ with $\dtw(G)=\Omega(\tw(G)\Delta(G))$.

A \textit{tree-partition} of a graph $G$ is a partition of $V(G)$ into (disjoint) bags $(B_i)_{i\in V(T)}$, where $T$ is a tree,  such that $uv\in V(G)$ implies that the bags of $u$ and $v$ are the same or adjacent in $T$. 
The width is the size of the largest bag, and the \textit{tree-partition-width} $\tpw(G)$ of $G$ is found by taking the minimum width over all tree-partitions of $G$.
The notion was first introduced by Seese~\cite{Seese85} under the name \emph{strong treewidth}.
\begin{observation}
\label{obs:related}
    For any graph $G$, $\tpw(G)-1\leq \dtw(G)\leq (\Delta+1) \tpw(G)$.
\end{observation}
The upper bound follows because after rooting the tree $T$ of a tree-partition of $G$, we may add a vertex $v$ to its parent bag $B_p$ (where $p$ is the parent of $c$ and $v\in B_c)$ whenever $v$ has a neighbour in $B_p$. This gives a domino tree decomposition of at most $\Delta+1$ times the width. 
For the lower bound, note that after rooting a domino tree decomposition, we can remove each vertex from a bag whenever it is not the bag closest to the root containing it. This yields a tree-partition, since for any edge $uv$, at least one of the two is still present in the highest bag they originally shared. (The construction was given by Bodlaender and Engelfriet~\cite{BodlaenderE97}.)

It was first shown that $\tpw(G)\leq 24(\tw(G)+1)\Delta(G)$ by an anonymous referee of \cite{DingO95} and the constant was later improved (to approximately $17.5$) by Wood~\cite{Wood09}. The lower bound $2\tpw(G)\geq \tw(G)+1$ was shown by Seese~\cite{Seese85}.
Tree-partitions for which the underlying tree has bounded maximum degree have been studied as well by Distel and Wood~\cite{Distel2024}.

Bodlaender, Groenland and Jacob \cite{BodlGroJacob} showed that \textsc{Domino Treewidth} with the domino treewidth as parameter is complete for the
parameterized complexity class XALP: deciding if a given graph has a tree decomposition of width at
most $k$ and spread at most $2$ is XALP-complete. They proved a similar result for \textsc{Tree-Partition-Width} and provided an FPT-approximation algorithm for this as well. 

The parameter \textit{treespan} can be defined in terms of the occupation time in a graph searching game~\cite{FominHT04} and is polynomially tied to domino treewidth (see \cite[Theorem 54]{Jacob25}).

Downey and McCartin~\cite{DowneyM04,DowneyM05} also studied the notion of spread for path decompositions under the name \textit{persistence}. They motivated this notion as being a more natural notion for a `path-like' graph, particularly for online computational problems. 
They showed amongst others that deciding if a given graph has a path decomposition of given width and persistence is W[$t$]-hard for all $t$ via a reduction from \textsc{Bandwidth}, with the sum of the width and persistence as parameter. This reduction can be easily turned into an XNLP-hardness proof, using the recent result
by Bodlaender, Groenland, Nederlof and Swennenhuis~\cite{BodlaenderGNS24} that
\textsc{Bandwidth} is XNLP-complete. An XNLP-membership proof can be done in the same fashion as the XNLP-membership of \textsc{Bandwidth} \cite{BodlaenderGNS24}. Thus
\textsc{Bounded Persistence Pathwidth} with sum of width
and persistence as parameter is XNLP-complete.

\section{Preliminaries}

We use the short-cut notation $[n]=\{1,\dots,n\}$.

The \emph{$(m\times n)$-grid} is the graph on vertex set $[m]\times [n]$ and edge set
\[
\{(i,j)(i+1,j): i\in[m-1],\ j\in[n]\}
\cup
\{(i,j)(i,j+1): i\in[m],\ j\in[n-1]\}.
\]

A {\em tree decomposition} of a graph $G$ is a pair $(T,(B_t)_{t\in V(T)})$ with
$T$ a tree and \textit{bags} $B_t\subseteq V(G)$ for each $t\in V(T)$  such that 
\begin{itemize}
    \item $\bigcup_{t\in V(T)} B_t = V(G)$,
    \item for all edges $vw\in E(G)$, there is an $t$ with $v,w\in B_t$, and 
    \item for all $v$, the nodes $\{t\in V(T)~|~v\in B_t \}$ form a subtree (denoted $T_v$) of $T$.
\end{itemize}
The {\em width} of the decomposition is $\max_{t\in V(T)} |B_t|-1$, and the {\em treewidth} of a graph $G$ is the minimum width over all tree decompositions of $G$. A
{\em path decomposition} is a tree decomposition $(T,(B_t)_{t\in V(T)})$ for which $T$ is a path and the {\em pathwidth} is the minimum width over all path decompositions of $G$. The \textit{spread} of a vertex $v$ in a
tree decomposition $(T,(B_t)_{t\in V(T)})$ is the number of bags $v$ is in: $|\{t \in V (T):v \in B_t\}|$. The \emph{spread} of a tree-decomposition $\mc{T} =(T, (B_t)_{t\in V(T)})$ is the maximum over $v \in V(G)$ of the number of bags of $\mc{T}$ containing $v$. We assume for simplicity that $B_t \neq \emptyset$ for each $t \in T$. It is easy to verify that we can do this without changing the width or the spread of a tree/path decomposition.

The following lemma is standard (see e.g. \cite[Corollary 9]{wood2025treedecompositionssmallwidth}). 
\begin{lemma}[Folklore]
\label{lem:standard}
    Let $(T,(B_t)_{t\in V(T)})$ be a tree decomposition of $G$. Let $w:V(G)\to \mathbb{R}_{\geq 0}$ be a weight function.
Then there is $t\in V(T)$ such that for each component $C$ of $G-B_t$, we find that $w(C)\leq w(V(G))/2$.
\end{lemma}

\section{Well-behaved spread using thrice the width}
In this section we prove Theorem~\ref{thm:spread_constant_min_width_upper_bound} which is immediately implied by the result below.
\begin{theorem}
\label{thm:spread_constant_min_width_upper_bound2}
Let $b\geq 5$ be an integer.
    Every graph $G$ admits a tree decomposition of width at most $(3+8/b)(\textup{tw}(G)+1)$ for which each vertex $v\in V(G)$ is in at most $2b(d(v)+1)$ bags. 
\end{theorem}
A key idea in the best bound of 14 from Wood~\cite{wood2025treedecompositionssmallwidth} is the definition of a \textit{slick} tree decomposition $(T,(B_x)_{x\in V(T)})$, in which the tree $T$ is rooted and the following property holds: if a vertex $v\in V(G)$ is in $B_x\cap B_y$ where $x$ is the parent of $y$, then $v$ has a neighbour in $B_y\setminus B_x$. This property automatically ensures that each vertex has spread at most  $d(v)+1$.

Using separators, a rooted tree decomposition can be constructed using induction on the number of vertices in the graph. 
After the first separator $S$ is obtained, a recursive call can be placed on the components of $G-S$; the aim is then to glue together the tree decompositions obtained via recursion using a bag containing $S$. At the next recursion depth, we are given a special set $S$ that needs to become part of the root bag (thinking ahead of our gluing operation), and we find a `balanced' separator $X$ for $S$, passing on $X$ and parts of $S$ to the next recursion depth as special set. Such an approximation scheme for treewidth has been used since the 1990s~\cite{BodlaenderGHK95,RobertsonS13}, see \cite{Bodlaender24} for a recent overview.

Wood~\cite{wood2025treedecompositionssmallwidth} cleverly noticed that in this routine, before passing a special set $S'$ to a recursive call, we can also replace it by a set $S$ where we add some neighbours of vertices in $S$ in order to ensure the tree decomposition will become slick. In order to keep the bag size smaller compared to their approach, we will only add neighbours for a `negligible' fraction of the vertices, while ensuring all vertices have a neighbour added once-in-a-while to upper bound the spread.

To facilitate the reasoning about the spread of the vertices, we introduce the notion of a \textit{$b$-marking} for a rooted tree decomposition $(T,(B_t))$ of $G$. This is a partial function $m:V(G)\times V(T)\to \{1,\dots,b\}$ such that:
\begin{itemize}
    \item $m(v,t)$ is defined if and only if $v\in B_t$;
    \item if $v$ is in $B_c\cap B_p$ for $c$ the child of $p$ in $T$, then either $m(v,c)<m(v,p)$ or $v$ has a neighbour in $B_c$ which it does not have in $B_p$;
    \item if $v\in B_\ell$ and $\ell$ does not have a child $c$ with $v\in B_c$ (and $\ell$ is not the root) then $v$ has a neighbour in $B_\ell$ that it does not have in the parent bag of $\ell$.  
\end{itemize}
Note that this last property can always be ensured, since if $T_v$ denotes the subtree of $T$ consisting of $t$ with $v\in B_t$, and $\ell$ is a leaf of $T_v$, then we may simply remove $v$ from $B_\ell$ if $v$ has no neighbours in $B_\ell$.

The use of this notion is via the following observation.
\begin{lemma}
\label{lem:bmarking_to_spread}
    If a rooted tree decomposition admits a $b$-marking, then each vertex $v$ is in at most $b(d(v)+1)$ bags.
\end{lemma}
\begin{proof}
Let $(T,(B_t)_{t\in V(T)})$ be a tree decomposition of $G$ with $b$-marking $m$. Let $v\in V(G)$. 
    Let $T_v$ be the subtree of $T$ corresponding to the  bags containing $v$. We call a node $t$ of $T_v$ \textit{happy} if there is a neighbour of $v$ in $B_t$ which was not in $B_p$ for $p$ the parent of $t$. (The root is also happy.) Note that:
    \begin{itemize}
        \item each leaf node of $T_v$ is happy (by the third marking property);
        \item any rooted path with $b$ vertices has at least one happy vertex.
    \end{itemize}
    The second property holds since such a path $t_1,\dots,t_b$ must pass from parent to child in each step and therefore the marking goes down by 1, or a happy node occurs, by the second marking property. 
    We use this to show that the fraction of nodes in $V(T_v)$ that is happy is at least $1/b$. Since the number of non-root happy nodes is at most the number of neighbours of $v$ (as $u$ is the reason a $t\in V(T_v)$ is happy only if $t$ is the closest vertex to the root with $u\in B_t$), this shows that $|V(T_v)|\leq b(d(v)+1)$.
    
    The remainder of the proof is straightforward. We initialise $T'=T_v$ and iteratively shrink it as follows. Starting from a leaf $v_1$ (which is happy by definition), we take a longest rooted path $v_k,v_{k-1},\dots,v_2,v_1$ such that $v_2,\dots,v_k$ are unhappy (allowing for $v_k=v_1$). In particular, the parent $v_{k+1}$ of $v_k$ (which could be the root) is happy and $v_k,\dots,v_2$ forms a rooted path of unhappy vertices, and so $k\leq b$. We change $T'$ by removing $v_1,\dots,v_k$ and attaching any neighbours of these vertices to $v_{k+1}$. Since $v_{k+1}$ is happy, this cannot create rooted paths on at least $b$ vertices without a happy vertex. Repeating this, we eventually end up with only the root, and reversing the process shows that the fraction of happy vertices is at least $1/b$.
\end{proof}
To prove Theorem~\ref{thm:spread_constant_min_width_upper_bound2}, we prove a slightly stronger claim by induction on $n$ below.
\begin{lemma}
\label{lem:create_marked_decomposition}
Let $b\geq 5$ and $k\geq 2$ be integers.
Let $G$ be a graph with treewidth at most $k-1$. Let $S\subseteq V(G)$ be given 
with
\[
2k+1\leq |S|\leq (2+8/b)k+1
\]together with $m':S\to \{1,\dots ,b+1\}$ such that at most $\lceil(1/b)(2+8/b)k\rceil$ vertices in $S$ are marked $i$ by $m'$ for all $i\in [b]$. 
Then $G$ admits a tree decomposition $(T,(B_x)_{x\in V(T)})$ such that
\begin{itemize}
    \item $|B_x|\leq (3+8/b)k+1$ for all $x\in V(T)$,
    \item $T$ is rooted in $r$ and $S\subseteq B_r$, 
    \item there is a $(b+1)$-marking $m:V(G)\times V(T)\to [b+1]$ for which $m(s,r)=m'(s)$ for all $s\in S$. 
\end{itemize} 
\end{lemma}
\begin{proof}
If $|V(G)|\leq (3+8/b)k+1$ then the lemma holds since we may put all vertices in the same bag.

Suppose that $|V(G)|> (3+8/b)k+1$ and that the lemma has been shown already for graphs on fewer vertices. Let $S,m'$ be given as in the lemma statement. We define a weight function $w$ on $G$ by giving all vertices in $S$ weight $1$ and the remaining vertices weight $0$. We apply Lemma~\ref{lem:standard} on $(G,w)$ using any tree decomposition for $G$ of optimal width, in order to find a vertex set $X$ with $|X|\leq \tw(G)+1\leq k$  such that each component $C$ of $G-X$ satisfies $|C\cap S|\leq |S|/2$. 

For each component $C$ of $G-X$, we add vertices to a set $S_C$ and define a marking $m_C':S_C\to [b+1]$ as follows:
\begin{itemize}
    \item For all vertices $v$ in $S\cap C$ or $X\cap S$ marked 1 by $m'$ that have a neighbour in $C\setminus S$, we add $v$ and a fixed neighbour $u$ of $v$ in $C\setminus S$ to $S_C$ and set $m'_C(v)=m'_C(u)=b+1$.
    \item All other vertices $s\in S\cap C$ or $X\cap S$ with a neighbour in $C\setminus S$ are added to $S_C$ with mark $m'_C(s)=m'(s)-1$. 
    \item All vertices $x\in X\setminus S$ with a neighbour in $C\setminus S$ are added to $S_C$ add with mark $m'_C(x)=b$.
\end{itemize}
In particular, vertices in $X$ or $S$ without neighbours in $C\setminus S$ are not added to $S_C$.
Note that 
\begin{align*}
|S_C|&\leq |S\cap C|+|X|+|\{v\in S:m'(v)=1\}|\\
&\leq \lfloor |S|/2\rfloor +k+\lceil(1/b)(2+8/b)k\rceil\\
&\leq 2k+\lfloor 1/2+4k/b\rfloor +\lceil 4k/b\rceil\\
&\leq (2+8/b)k+1.
\end{align*}
where we used that $8/b\leq 2$ for $4\leq b$ and that $\lfloor 1/2+x\rfloor +\lceil x\rceil\leq 2x+1$ for all  real $x\geq0$.

We wish to apply the induction hypothesis on $G_C=G[(C\setminus S)\cup S_C]$ with special set $S_C$ and the marking $m'_C$ constructed above. If $|V(G_C)|\leq 3k$, then this is not needed and we may simply define $R_C=S_C$. Otherwise, we can ensure that $|S_C|\geq 2k+1$; if needed, we can add vertices from $C$ to $S_C$ and mark these $b+1$ by $m_C'$. Note that $|V(G_C)|<|V(G)|$ since at most $|S\cap C|+|X\cap S|\leq |S|/2+k< |S|$
vertices of $S$ remain present in $G_C$.

In order to ensure the marking $m'_C$ for $S_C$ satisfies the assumptions in the lemma statement, we can decrease some of the marks if needed. Any vertex in $S_C$ marked $i$ by $m_C'$ for some $i\in [b-1]$, is marked $i+1$ in $S$ by $m'$, so there are not too many of those. Moreover, at most \[
|S\cap C|+ |X|\leq \lfloor |S|/2\rfloor +k\leq 2k+\lfloor 1/2+4k/b\rfloor\leq 2k+8k/b
\]
vertices can ever receive a mark $\leq b$.
This means that if too many vertices received mark $b$ at this point, we can safely decrease some of their marks (`increasing their priority'). 

We are now ready to apply the induction hypothesis on $G_C=G[(C\setminus S)\cup S_C]$ with special set $S_C$ and the marking $m'_C$ adjusted as above (if needed). 
This yields a tree decomposition $T_C$ with root bag $R_C$ containing $S_C$ for each component $C$ together with a $(b+1)$-marking $m_C$. 

We make one additional bag $B_r=S\cup X$ and make it adjacent to the various $R_C$ bags. Note that $|B_r|\leq |S|+k\leq (3+8/b)k+1$.
This creates a tree decomposition which we root in $B_r$. We let the marking $m$ be defined by
\begin{itemize}
    \item $m(s,r)=m'(s)$ for all $s\in S$;
    \item $m(x,r)=b+1$ for all $x\in X\setminus S$;
    \item $m(v,t)=m_C(v,t)$ for $t\in V(T_C)$ and $v\in B_t$. 
\end{itemize}
In order to show this is a $(b+1)$-marking, we only need to consider what happens to $v\in B_p\cap B_c$ where $B_p=B_r=S\cup X$ is the root bag and $B_c=R_C$ for some component $C$. If $v\in X\setminus S$, then
\[
m(v,c)=m_C(v,c)=m'_C(v)\leq b=(b+1)-1=m(v,p)-1
\]
as desired. If $v\in S$ then 
\[
m(v,c)=m_C(v,c)=m'_C(v)= m'(v)-1=m(v,p)-1,
\]
unless $m'(v)=1$ in which case a neighbour of $v$ in $C$ which was not in $B_r$ will have been added. This shows that $m$ is a $(b+1)$-marking, as desired.
\end{proof}
\begin{proof}[Proof of Theorem~\ref{thm:spread_constant_min_width_upper_bound2}] 
As the case $\textup{tw}(G)=0$ is immediate, we assume  $\textup{tw}(G) \geq 1$.
    Let $b\geq 5$ be an integer and let $G$ be a graph and let $k=\textup{tw}(G)+1\geq 2$. The statement is immediately true when $|V(G)|\leq 2k+1$. Otherwise, we let $S$ contain $2k+1$ vertices of $G$ and set $m'(s)=b+1$ for all $s\in S$.  Lemma~\ref{lem:create_marked_decomposition} now provides a tree decomposition of width at most $(3+8/b)k=(3+8/b)(\text{tw}(G)+1)$ that admits a $(b+1)$-marking. By Lemma~\ref{lem:bmarking_to_spread} each vertex $v$ has spread at most $(b+1)(d(v)+1)\leq 2b(d(v)+1)$ in this. 
 \end{proof}

\section{Grid graphs and extensions}

We first provide a near-optimal width tree decomposition for grids with constant spread. 
\begin{lemma}
\label{lem:good_spread_width_for_grids}
Let $c>1$ and let $m,n\in \mathbb{N}$. Then the $(m\times n)$-grid has a path decomposition of width $cn$ and spread at most $\lceil1/(c-1)\rceil +1$.
\end{lemma}
We provide a `proof by picture' below and include a formal proof in the appendix.
\begin{center}
\begin{tikzpicture}[scale=0.55]

% parameters for the figure
\def\n{3}      % number of rows  (n)
\def\m{5}       % number of cols  (m)
\def\ai{1}      % start of A_i
\def\aii{2}     % start of A_{i+1}
\def\aiii{3}     % start of A_{i+1}
\def\b{\n}       % end of A_{[i,b]}

% Draw background grid (light gray)
\foreach \i in {1,...,\n}{
  \foreach \j in {1,...,\m}{
    \draw (\j,-\i) rectangle ++(1,1);
  }
}

%%%%%%%%%%%%%%%%%%%%%%%%%%%%%%%%%%%%%%%%%%%%%%%%%%%%%%%%%%%%
% Bag B_{i,j} at (j, j+1) rows (horizontal axis is i)
%%%%%%%%%%%%%%%%%%%%%%%%%%%%%%%%%%%%%%%%%%%%%%%%%%%%%%%%%%%%

% Choose a specific column j0 and j0+1 for illustration
\def\jone{1}   % j
\def\jtwo{2}   % j+1

% B_{i,j}:   A_[i,b] in column j
\foreach \i in {\ai,...,\b}{
  \fill[blue!35] (\jone,-\i) rectangle ++(1,1);
}

% B_{i,j}:   A_[1,i] in column j+1 (Flipped vertically)
\foreach \i in {1,...,\ai}{
  \pgfmathtruncatemacro{\flipped}{\n-\i+1}
  \fill[red!35] (\jtwo,-\flipped) rectangle ++(1,1);
}

%%%%%%%%%%%%%%%%%%%%%%%%%%%%%%%%%%%%%%%%%%%%%%%%%%%%%%%%%%%%
% Next bag B_{i+1,j} shifted right for clarity
%%%%%%%%%%%%%%%%%%%%%%%%%%%%%%%%%%%%%%%%%%%%%%%%%%%%%%%%%%%%

\begin{scope}[xshift=6cm]

% redraw grid
\foreach \i in {1,...,\n}{
  \foreach \j in {1,...,\m}{
    \draw (\j,-\i) rectangle ++(1,1);
  }
}

% B_{i+1,j}: A_[i+1,b] in column j
\foreach \i in {\aii,...,\b}{
  \fill[blue!35] (\jone,-\i) rectangle ++(1,1);
}

% B_{i+1,j}: A_[1,i+1] in column j+1 (Flipped vertically)
\foreach \i in {1,...,\aii}{
  \pgfmathtruncatemacro{\flipped}{\n-\i+1}
  \fill[red!35] (\jtwo,-\flipped) rectangle ++(1,1);
}

\end{scope}

%%%%%%%%%%%%%%%%%%%%%%%%%%%%%%%%%%%%%%%%%%%%%%%%%%%%%%%%%%%%
% Next bag B_{i+2,j} shifted right for clarity
%%%%%%%%%%%%%%%%%%%%%%%%%%%%%%%%%%%%%%%%%%%%%%%%%%%%%%%%%%%%

\begin{scope}[xshift=12cm]

% redraw grid
\foreach \i in {1,...,\n}{
  \foreach \j in {1,...,\m}{
    \draw (\j,-\i) rectangle ++(1,1);
  }
}

% B_{i+1,j}: A_[i+1,b] in column j
\foreach \i in {\aiii,...,\b}{
  \fill[blue!35] (\jone,-\i) rectangle ++(1,1);
}

% B_{i+1,j}: A_[1,i+1] in column j+1 (Flipped vertically)
\foreach \i in {1,...,\aiii}{
  \pgfmathtruncatemacro{\flipped}{\n-\i+1}
  \fill[red!35] (\jtwo,-\flipped) rectangle ++(1,1);
}

\end{scope}

%%%%%%%%%%%%%%%%%%%%%%%%%%%%%%%%%%%%%%%%%%%%%%%%%%%%%%%%%%%%
% Next bag B_{i+2,j} shifted right for clarity
%%%%%%%%%%%%%%%%%%%%%%%%%%%%%%%%%%%%%%%%%%%%%%%%%%%%%%%%%%%%

\begin{scope}[xshift=18cm]

% redraw grid
\foreach \i in {1,...,\n}{
  \foreach \j in {1,...,\m}{
    \draw (\j,-\i) rectangle ++(1,1);
  }
}

% B_{i+1,j}: A_[i+1,b] in column j
\foreach \i in {\ai,...,\b}{
  \fill[blue!35] (\jone+1,-\i) rectangle ++(1,1);
}

% B_{i+1,j}: A_[1,i+1] in column j+1 (Flipped vertically)
\foreach \i in {1,...,\ai}{
  \pgfmathtruncatemacro{\flipped}{\n-\i+1}
  \fill[red!35] (\jtwo+1,-\flipped) rectangle ++(1,1);
}

\end{scope}

\end{tikzpicture}
\end{center}
The picture above shows four consecutive bags for $c=1+1/3$ via the coloured regions. Each column denotes a single column of the $(m\times n)$-grid, whereas any row in the picture corresponds to the union of at most $\lceil n/3\rceil $ consecutive rows of the grid.

We now move to the proof of \cref{thm:spread_constant_min_width_lower_bound}. We will prove something stronger.

\begin{theorem}\label{lowermain} For all positive integers $s,w_0 > 0$, there exists a graph $G$ of treewidth $w \geq w_0$ such that $\Delta(G) = 3$ and in every tree-decomposition with width at most $3w-8$, some vertex $v \in V(G)$ has spread strictly more than $s$.\end{theorem}

Note that \cref{lowermain} not only implies \cref{thm:spread_constant_min_width_lower_bound}, but also reveals the following: for a graph with treewidth $w$, if we only consider tree-decompositions of width less than $3w - o(w)$, we cannot  guarantee existence of a tree-decomposition with spread bounded by \textit{any} function of the degree, let alone some linear function.

\begin{proof}[Proof of \cref{lowermain}.]
For simplicity, we prove  a slight weakening of \cref{lowermain} with the condition $\Delta(G) = 3$ replaced by $\Delta(G) = 4$. To recover \cref{lowermain} as stated, one needs to replace the grid graph in our example by a ``wall''. The modification is routine, but complicates the notation significantly, thus we leave the details to the motivated reader.

Given positive integers $w_0$ and $s$, let $n=w_0$ and consider a modification of the $((m^2+1) \times (n+1))$-grid, with $m$ to be chosen later. Formally, consider $G$ where \begin{align*}V(G) &= [m^2+1] \times [n+1] \cup [m + 1] \qquad \mathrm{and} \\
E(G) &= \set{(i, j)(i+1,j) : i \in [m^2]} \cup \set{(i, j)(i, j+1) : j \in [n]} \\ &\cup \set{(1,n+1)(m^2+1,n+1)} \\ &\cup \set{(i)(i+1) : i \in [m]} \\ &\cup \set{(i) (i m+1, n+1) : i \in [m]}.\end{align*} Essentially, to obtain this graph, we first add an edge between the top two corners to the $((m^2+1) \times (n+1))$-grid. Then, we add an extra path above the grid with only $m + 1$ vertices, with edges to $m + 1$ vertices in the original grid that are spaced $m$ vertices apart. See Figure~\ref{f:example} for a visualization of $G$. Clearly, $\Delta(G) = 4$.

\begin{figure}[ht!]
\centering
\begin{tikzpicture}[
  vertex/.style={circle, draw, inner sep=1.5pt},
  gridedge/.style={draw, line width=0.5pt},
  extraedge/.style={draw, line width=1pt}
]

\def\ncols{10}
\def\nrows{3}
\def\dx{1.0}
\def\dy{1.0}

\foreach \r in {1,...,\nrows} {
  \pgfmathtruncatemacro{\y}{\nrows-\r} % 2,1,0
  \foreach \c in {1,...,\ncols} {
    \node[vertex] (g-\r-\c) at ({(\c-1)*\dx},{\y*\dy}) {};
  }
}

\foreach \r in {1,...,\nrows} {
  \foreach \c in {1,...,\numexpr\ncols-1\relax} {
    \draw[gridedge] (g-\r-\c) -- (g-\r-\the\numexpr\c+1\relax);
  }
}

\foreach \r in {1,...,\numexpr\nrows-1\relax} {
  \foreach \c in {1,...,\ncols} {
    \draw[gridedge] (g-\r-\c) -- (g-\the\numexpr\r+1\relax-\c);
  }
}

\draw[extraedge, bend left=50] (g-1-1) to (g-1-\ncols);

\foreach \i/\col in {1/1, 2/4, 3/7, 4/10} {
  \coordinate (topcol-\i) at ($(g-1-\col) + (0,\dy)$);
  \node[vertex] (p-\i) at (topcol-\i) {};
}

\foreach \i in {1,...,3} {
  \draw[extraedge] (p-\i) -- (p-\the\numexpr\i+1\relax);
}

\foreach \i/\col in {1/1, 2/4, 3/7, 4/10} {
  \draw[extraedge] (p-\i) -- (g-1-\col);
}

\end{tikzpicture}
\caption{The graph $G$ described above is depicted for $m = 3$ and $n = 2$}\label{f:example}
\end{figure}

\begin{claim}\label{tw}
$\mathrm{tw}(G) \leq n+3$
\end{claim}

\begin{proof}
First, there is an easy argument that the $(m \times n)$-grid has treewidth at most $n$: take bags $B_{i,j} \coloneq \set{(i, j), (i, j+1), ..., (i, n), (i+1, 1), (i+1, 2), ... (i+1, j)}$ for $1 \leq i \leq m-1$ and $1 \leq j \leq n$. To form a tree-decomposition of width $n$, connect $B_{i,j}$ to $B_{i, j+1}$ for all $1 \leq i \leq m-1$ and $1 \leq j \leq n-1$, then connect $B_{i,n}$ to $B_{i+1, 1}$ for all $1 \leq i \leq m-2$. Let $H$ be the graph obtained from $G$ by removing the additional edge $(1,n+1)(m^2+1,n+1)$ that was placed on the $(m^2+1)\times (n+1)$-grid. As $H$ is a subgraph of the
$(m^2+1)\times (n+1)$-grid, we have $\mathrm{tw}(H) \leq n+2$. As $G$ is obtained from $H$ by adding one edge, the claim holds.
\end{proof}

\begin{claim}\label[claim]{diam}
$\mathrm{diam}(G) \leq 3m + 2n + 2$
\end{claim}
\begin{proof}
To construct a path from vertex $u$ to vertex $v$ in $G$, we can follow at most $n$ edges from $u$ up to the top of the main $((m^2+1) \times (n+1))$-grid, at most $m$ edges to a vertex that connects to the top path, then 1 edge to reach the top path (if $u$ is not already in the top path). Then, either it takes at most $m$ edges to reach $v$ if $v$ is in the top path, or it takes at most $m$ edges to reach a nearest column to $v$ with a vertex connecting to the top path. Finally, if $v$ is not in the top path, it takes 1 edge to get back down to the main grid, at most $m$ edges to reach the column with $v$, and at most $n$ more edges to reach the correct row.
\end{proof}

Suppose for a contradiction that there is a tree-decomposition $\mc{T} = (T, \set{B_t}_{t \in V(T)})$ of $G$ with width at most $3n+1$ and with spread at most $s$.

\begin{claim}\label[claim]{pathlen}
Let $P$ be any path in $T$. Then $|V(P)| \leq (3m+2n+3)s$.
\end{claim}

\begin{proof}
Consider any path $P$ in $T$ from $t_1 \in V(T)$ to $t_2 \in V(T)$. Pick any vertices $u \in B_{t_1}$ and $v \in B_{t_2}$. Let $P_{uv}$ a shortest path between $u$ and $v$ in $G$. By \cref{diam}, $P_{uv}$ has at most $3m+2n+3$ vertices. Further, each bag of a vertex of $P$ must contain a vertex of $P_{uv}$ as a standard property of tree-decompositions. Since each vertex of $P_{uv}$ can appear in at most $s$ bags, $|V(P)| \leq (3m + 2n + 3)s$.
\end{proof}

For any connected subgraph $J$ of $G$, let $T_J$ be the subtree $T_J \subseteq T$ such that $V(T_J) = \set{t \in V(T) : B_t \cap J \neq \emptyset}$. Consider the subgraph $H \coloneq G[\set{(1,j) : 1 \leq j \leq n+1} \cup \set{(m^2+1,j) : 1 \leq j \leq n+1} \subseteq G$, i.e. the graph including the leftmost and rightmost columns of the main grid and the edge between the top two corners. As $|V(H)| = 2n+2$, at most $(2n+2)s(3n+2)$ vertices of $G$ can appear in bags of $T_H$. 
Let $k=(3n+2)(3m+2n+3)s+2$. Assume \begin{equation}\label{e:mLower}
m^2-1 \geq k((2n+2)s(3n+2)+1).\end{equation} Let $C_i=G[\set{(i, j) : 1 \leq j \leq n}]$, representing the connected graph induced on column $i$ of $G$ (not including the top path of length $m$).  Let $C_{i,j} = G[\bigcup_{i \leq \ell \leq j} \set{(\ell, r) : 1 \leq r \leq n}]$, representing a subgrid of the $((m^2+1) \times (n+1))$-grid used to construct $G$. By \eqref{e:mLower} there must be some $2 \leq i < j \leq m^2$ such that the subtree $T_{C_{i,j}}$ is disjoint from $T_H$ and $j-i+1 \geq k$, i.e. $C_{i,j}$ is a subgrid that consists of at least $k$ columns and has no vertices in the same bag as any vertex of $H$.

Clearly $T_{C_r} \subseteq T_{C_{i,j}}$ for any $i \leq r \leq j$, so $T_{C_r}$ is disjoint from $T_H$ for any $i \leq r \leq j$. Consider in particular $T_{C_\ell}$ such that $\ell = \arg \max_{i \leq r \leq j} d_T(T_H, T_r)$. Consider the shortest path $P$ from $T_H$ to $T_{C_\ell}$. By \cref{pathlen}, $|V(P)| \leq (3m + 2n + 3)s$, and each bag of $T$ contains at most $3n+2$ vertices of $G$, so at most $(3n+2)(3m + 2n + 3)s = k-2$ choices of $r$ subject to $i \leq r \leq j, r \neq \ell$ have $T_{C_r} \cap P \neq \emptyset$. Since $(j-i+1) \geq k$, there must be some  $r \neq \ell$ with $i \leq r \leq j$ such that $T_{C_r} \cap P = \emptyset$. Fix this $r$.  

Then, $T_{C_r}$ cannot intersect $T_H$ or the shortest path in $T$ from $T_H$ to $T_{C_\ell}$, $T_{C_\ell}$ cannot intersect  the shortest path from $T_H$ to $T_{C_r}$ by the choice of $\ell$ and $r$, and $T_H$ cannot intersect the shortest path from $T_{C_r}$ to $T_{C_\ell}$ since $T_H$ and $T_{C_{i,j}}$ are disjoint, $T_{C_{i,j}}$ is connected  and $T_{C_r},T_{C_\ell} \subseteq T_{C_{i,j}}$.  Therefore, there must be some vertex $t' \in T$ such that $t'$ intersects each of the shortest path from $T_H$ to $T_{C_\ell}$, the shortest path from $T_H$ to $T_{C_r}$, and the shortest path from $T_{C_\ell}$ to $T_{C_r}$, with $t'$ not contained in any of these three subtrees.

Suppose without loss of generality that $r < \ell$. As $B_{t'}$ must contain one vertex of each path from $H$ to $C_r$,   $B_{t'}$ must contain a vertex of each of the rows of the main grid with the second coordinate strictly between $1$ and $r$. Similarly, by considering paths from $C_r$ to $C_\ell$, we conclude that $B_{t'}$ must contain one vertex of a vertex of each of the rows of the main grid from columns $r+1$ to $\ell-1$. Finally, by considering paths from $C_\ell$ to $H$,  we get further $n+1$ vertices in $B_{t'}$ from columns $\ell+1$ to $m^2$. Thus $|B_{t'}| \geq 3n+3$ contradicting the choice of $\mc{T}$.

It remains to show that $m$ can be chosen to satisfy \eqref{e:mLower}. We can choose $m$ large enough such that $\frac15 \sqrt m >((2n+2)s(3n+2)+1)$, \\ $\sqrt m > (3n+2)s$, $m > 2n + 3$. Then, we only need \\ $m^2 - 1 \geq \frac15 \sqrt m (\sqrt m(4m) + 2)$, which is easily seen to be true for large enough $m$.

By \cref{tw}, we thus find a graph $G$ of treewidth  $n+1 \leq w \leq n+3$ such that $\Delta(G) = 4$ and in every tree-decomposition with width at most $3n+1 \geq 3w-8$, some vertex has spread strictly more than $s$.
\end{proof}

\section{Proof of \cref{thm:grid_linear_spread}}

Let $G$ be the $(n \times n)$-grid for this section. We denote by $L,R,F$ and $C$  the \emph{sides} of the grid, i.e.
\begin{align*}
    L&= \{1\} \times [n],\\ 
    R&= \{n\} \times [n], \\
    F&= [n] \times \{1\}, \\
    C&= [n] \times \{n\}.
\end{align*}

We say that a subgraph $H$ of $G$ is \emph{global} if $V(H)$ intersects every side of $G$, and \emph{non-global}, otherwise.	We say that a subset $X$ of $V(G)$ is \emph{connected} (resp. global, non-global)  if $G[X]$ is. Let $\mc{T}=(T, \{B_t\}_{t \in V(T)})$ be a tree-decomposition of $G$, let $n + a - 1$ be the width of $\mc{T}$, i.e. $n+a = \max_{t \in V(T)} |B_t|$. For an oriented edge $tt' \in E(T)$, let $T_{tt'}$ denote the component of $T \setminus tt'$ containing $t'$, and let $V_{tt'} = \cup_{v \in V(T_{tt'})} B_v$ be the set of vertices of $G$ which belong to the bags in this component.

     \begin{obs}\label[obs]{obs:conn} Given a tree-decomposition $\mc{T}$ of a graph $G$, there exists a tree-decomposition $\mc{T}' = (T', \{B_t\}_{t \in V(T')})$ of $G$ with the same width and spread such that for any $tt' \in E(T')$, $V_{tt'}$ is connected. \end{obs}

     \begin{proof}
     Given $tt' \in E(T)$, let $H_1, H_2, ..., H_k$ be the components of $V_{tt'}$. We can construct $T'$ by replacing $T_{tt'}$ with $k$ copies of $T_{tt'}$, say $T_{tt'1}, ..., T_{tt'k}$. Similarly, let $t'_1, ..., t'_k$ be the $k$ copies of $t'$. For each copy $i$, add an edge between $t'_i$ and $t$. Then, for a vertex $x \in T_{tt'i}$, define $B_x = B_{T_{tt'}} \cap V(H_i)$. It is routine to verify that this construction satisfies all necessary properties of a tree decomposition, and for any vertex $t'_i$, $V_{tt'_i}$ is connected. One can simply repeat this process as needed to obtain the full construction of $\mc{T'}$.
     \end{proof}
		
	 By \cref{obs:conn}, we may assume that $V_{tt'}$ is connected for every $tt' \in E(T)$. Note that global connected sets are pairwise intersecting, which implies the following observation.
	 
	 \begin{obs}\label[obs]{obs:separator} There exists $t \in V(T)$ such that $V_{tt'} - B_t$ is non-global for every $tt' \in E(T)$. \end{obs}	
	 
	 \begin{figure}[ht!]
	 	\centering
	 	 \begin{tikzpicture}[scale=0.6]
	 		% styles
	 		\tikzset{
	 			ord/.style={circle,draw,fill=white,inner sep=1.5pt}, % ordinary (white) vertex
	 			Xv/.style={circle,draw,fill=black,inner sep=2pt}  % vertex in X (black)
	 		}
	 		
	 		% create nodes (1..8)
	 		\foreach \i in {1,...,8}{
	 			\foreach \j in {1,...,8}{
	 				\node (v\i\j) at (\i,\j) [ord] {};
	 			}
	 		}
	 		
	 		% horizontal edges
	 		\foreach \i in {1,...,7}{
	 			\foreach \j in {1,...,8}{
	 				\pgfmathtruncatemacro{\ip}{\i+1}
	 				\draw (v\i\j) -- (v\ip\j);
	 			}
	 		}
	 		
	 		% vertical edges
	 		\foreach \i in {1,...,8}{
	 			\foreach \j in {1,...,7}{
	 				\pgfmathtruncatemacro{\jp}{\j+1}
	 				\draw (v\i\j) -- (v\i\jp);
	 			}
	 		}
	 		
	 		% Fill X: first part { (i,j) | j in {2,3,4} } for all i
	 		\foreach \j in {2,3,5}{
	 			\foreach \i in {1,...,8}{
	 				\node at (v\i\j) [Xv] {};
	 			}
	 		}
	 		
	 		% Fill X: second part { (i,j) | i in {2,3,5}, j <= 5 }
	 		\foreach \i in {2,3,5}{
	 			\foreach \j in {1,...,5}{
	 				\node at (v\i\j) [Xv] {};
	 			}
	 		}
	 	\end{tikzpicture}
	 	\hspace{1cm}
	 		\begin{tikzpicture}[scale=0.6]
	 			% style for ordinary (white) vertices
	 			\tikzset{ord/.style={circle,draw,fill=white,inner sep=1.5pt}}
	 			% style for highlighted vertices (in X)
	 			\tikzset{Xv/.style={circle,draw,fill=black,inner sep=2pt}}
	 			
	 			% create and name nodes (1..8)
	 			\foreach \x in {1,...,8}{
	 				\foreach \y in {1,...,8}{
	 					\node (v\x\y) at (\x,\y) [ord] {};
	 				}
	 			}
	 			
	 			% draw horizontal edges
	 			\foreach \x in {1,...,7}{
	 				\foreach \y in {1,...,8}{
	 					\pgfmathtruncatemacro{\xp}{\x+1}
	 					\draw (v\x\y) -- (v\xp\y);
	 				}
	 			}
	 			
	 			% draw vertical edges
	 			\foreach \x in {1,...,8}{
	 				\foreach \y in {1,...,7}{
	 					\pgfmathtruncatemacro{\yp}{\y+1}
	 					\draw (v\x\y) -- (v\x\yp);
	 				}
	 			}
	 			
	 			% Definition of I
	 			% For each i in I, fill (j,i) for j<i and (i,j) for j>i
	 			\foreach \i in {3,5}{
	 				% (j,i) with j < i
	 				\foreach \j in {1,...,7}{
	 					\pgfmathtruncatemacro{\cond}{ifthenelse(\j < \i,1,0)}
	 					\ifnum\cond=1
	 					% place a node on top of the named node so it becomes black-filled
	 					\node at (v\j\i) [Xv] {};
	 					\fi
	 				}}
 				\foreach \i in {2,3}{
	 				% (i,j) with i < j
	 				\foreach \j in {2,...,8}{
	 					\pgfmathtruncatemacro{\cond}{ifthenelse(\i<\j,1,0)}
	 					\ifnum\cond=1
	 					\node at (v\i\j) [Xv] {};
	 					\fi
	 				}
	 			}
	 	\end{tikzpicture}
	 	\caption{A straight (left) and a crooked (right) $(3 \times 3)$-semigrid in an $(8 \times 8)$-grid defined by the sets $I=J = \{2,3,5\}$}
	 	\label{fig:tikz1}
	 \end{figure}

	 We now introduce the key objects in our argument, which we call semigrids.
	 For a non-empty, finite set $J$ of integers let $\max J$ and $\min J$ denote its maximum and minimum elements, respectively.
	 
	  Let $1 \leq k \leq n$. We say that $X \subseteq  V(G)$ is a \emph{straight $(k \times k)$-semigrid} if there exist $I, J \subseteq [n]$ with $|I|=|J| =k $ such that $$X = \{(i,j) \: | \: i \in [n]  j \in J\} \cup \{(i,j) \: | i \in I, j \leq \max J \:\}.$$
	We also call $X$ a straight  $(k \times k)$-semigrid if it is obtained from the set defined above by an automorphism of $G$, i.e. by ``rotating the grid".
	
	We say that $X \subseteq V(G)$ is a \emph{crooked $(k \times k)$-semigrid} if there exists $I \subseteq [n]$ with $|I|=k$ such that $$X = \{(j,i) \: | \: j < i, i \in I- \{\min I\} \} \cup \{ (i,j) \: | \: i \in I - \{\max I\}, i < j \}.$$ 
	Similarly,  $X$ is also a crooked  $(k \times k)$-semigrid if it is obtained from the set defined above by a rotation.
	
	A \emph{$(k \times k)$-semigrid} is a straight or crooked $(k \times k)$-semigrid. Note that every $(k \times k)$-semigrid is connected.  
	
	For $U \subseteq V(G)$ let $\bd{U}$ denote the set of all vertices $v \in U$ with a neighbour in $V(G)-U$. We say that a pair $(U,Z)$ of subsets of $V(G)$ is \emph{$k$-good} if \begin{itemize}
		\item $\bd{U} \subseteq Z \subseteq U$,
		\item $U$  is connected,
		\item $U-Z$ is non-global, and 
		\item $U-Z$ contains a $(k \times k)$-semigrid. 
	\end{itemize}
	
	\begin{lemma}\label{lem:l_boundary}
		Let $(U,Z)$ be $k$-good for some $k \geq 3$. Then $|Z| \geq n$.
	\end{lemma}
	
	\begin{proof}
		Let $X$ be a $(k \times k)$-semigrid in $U-Z$. We may assume by symmetry that $X \cap L, X \cap C \neq \emptyset $ and that $U-Z \cap F = \emptyset$, as $U-Z$ is non-global. If $U \cap R \neq \emptyset$ then $U$ contains a vertex in every column of $G$, and so $Z$ contains a vertex in every column and the lemma holds.
		
		Thus we may assume that $U \cap R = \emptyset$. Thus  $X$ is crooked and as $k \geq 3$ there exists $i_0 \in [n]$ such that 
		\begin{equation}\label{e:crooked} \{(j,i_0) \: | \: j < i_0 \} \cup \{(i_0, j) \: | \:   j > i_0 \} \subseteq X. \end{equation}
		 Let \begin{align*} A_j = \begin{cases} \{ (j,i) \: | \: i \leq i_0\}  & \mathrm{for} \qquad j \leq i_0, \\  \{(i,j) \: | \: i \geq i_0\} & \mathrm{for} \qquad i_0 < j \leq n. \end{cases} \end{align*} (See Figure~\ref{fig:tikz2}.)
		  Note that the sets $A_1,\ldots,A_n$ are connected and pairwise vertex disjoint. Each of them contains a vertex in $V(G) - (U-Z)$, as $F \cup R \subseteq V(G)-(U-Z)$ by our assumptions. Moreover, each of these sets contains a vertex of $U$ as each of them, except possibly $A_{i_0}$, contains a vertex of $X$ (and thus of $U$) by \eqref{e:crooked}. Finally, $(i_0,i_0) \in A_{i_0} \cap U$  as $(i_0,i_0)$ is adjacent to a vertex of $X \subseteq U-Z$ and $\bd{U} \subseteq Z$.  It follows that $A_i \cap Z \neq \emptyset$ for every $i \in [n]$, again implying the lemma.
	\end{proof}
 
	\begin{figure}[ht!]
		\centering
		\begin{tikzpicture}[scale=0.9]
			% --- parameters ---
			\def\pad{0.35} % padding for the orange boxes
			
			% --- styles ---
			\tikzset{
				ord/.style={circle,draw,fill=white,inner sep=1.5pt},
				blk/.style={circle,draw,fill=black,inner sep=2.0pt},
				Abox/.style={draw=orange, line width=0.35pt, fill=orange, fill opacity=0.22, draw opacity=0.75}
			}
			
			% --- orange boxes for A_j + labels ---
			% A_j = {(j,i) : i<=3} for j=1,2,3
			\foreach \j in {1,2,3}{
				\path[Abox] (\j-\pad, 1-\pad) rectangle (\j+\pad, 3+\pad);
				% label BELOW the box, in black
				\node[black, anchor=north] at (\j, 1-\pad) {$A_{\j}$};
			}
			
			% A_j = {(i,j) : i>=3} for j=4,5,6,7,8
			\foreach \j in {4,5,6,7,8}{
				\path[Abox] (3-\pad, \j-\pad) rectangle (8+\pad, \j+\pad);
				% label to the RIGHT of the box, in black
				\node[black, anchor=west] at (8+\pad, \j) {$A_{\j}$};
			}
			
			% --- nodes (1..8) ---
			\foreach \x in {1,...,8}{
				\foreach \y in {1,...,8}{
					\node (v\x\y) at (\x,\y) [ord] {};
				}
			}
			
			% --- edges ---
			\foreach \x in {1,...,7}{
				\foreach \y in {1,...,8}{
					\pgfmathtruncatemacro{\xp}{\x+1}
					\draw (v\x\y) -- (v\xp\y);
				}
			}
			\foreach \x in {1,...,8}{
				\foreach \y in {1,...,7}{
					\pgfmathtruncatemacro{\yp}{\y+1}
					\draw (v\x\y) -- (v\x\yp);
				}
			}
			
			% --- black vertices: {(j,3) : j<=3} U {(3,j) : j>3} ---
			\foreach \j in {1,2}{
				\node at (v\j3) [blk] {};
			}
			\foreach \j in {4,5,6,7,8}{
				\node at (v3\j) [blk] {};
			}
		\end{tikzpicture}
		
		\caption{Sets $A_j$ in the proof of \cref{lem:l_boundary}.}
		\label{fig:tikz2}
	\end{figure}
	
	We say that an oriented edge $tt' \in E(T)$ is \emph{$k$-good} if $(V_{tt'}, B_t \cap B_{t'})$ is $k$-good.
	
	\begin{lemma}\label{lem:l_step} Let $k \geq  2a+2$ be an integer and let $tt' \in E(T)$ be $k$-good. Then there exists $t'' \in T_{tt'}$ such that $t't'' \in E(T)$ is $(k-2a)$-good. 
	\end{lemma}
	
	\begin{proof} Let $B' = B_{t'} - B_t$.
		By \cref{lem:l_boundary}, $|B_t \cap B_{t'}| \geq n$. Thus  $|B'| \leq a$. Let $X$ be a  $(k \times k)$-semigrid  with $X \subseteq V_{tt'} - B_t$. Let $k' = k-2a$. It is easy to see that $X$ contains a $(k' \times k')$-semigrid $X'$ such that $X' \subseteq X - B' \subseteq V_{tt'} - B_{t'}$. As $X'$ is connected there exists $t'' \in V(T_{tt'})$ such that  $t't'' \in E(T)$  and $X' \subseteq V_{t't''} - B_{t'}$.  To show that $t't''$ is $k'$-good it remains to verify that
		\begin{itemize}
			\item $\bd{V_{t't''}} \subseteq  B_{t''} \cap B_{t'} \subseteq V_{tt''}$,
			\item $V_{t't''}$  is connected, and
			\item $V_{t't''}-(B_{t''} \cap B_{t'})$ is non-global. 
		\end{itemize}
		The first of these conditions holds by definition of a tree decomposition, and the second by our assumption on the tree-decomposition. Finally, the third condition holds as $$V_{t't''}-(B_{t''} \cap B_{t'}) = V_{t't''} - B_{t'} \subseteq V_{tt'}-B_{t'} \subseteq V_{tt'}-(B_{t'}\cap B_t),$$
		where the last set is non-global as  $(V_{tt'}, B_t \cap B_{t'})$ is $k$-good.
	\end{proof}
	
	\begin{cor}\label{c:induct} Let $s \geq 0$  be an integers, such that $sa \leq n$  and let $tt' \in E(T)$ be $(2sa+3)$-good. Then the spread of $\mc{T}$ is at least $\frac{s}{2}$.
	\end{cor}
	\begin{proof} By repeatedly applying \cref{lem:l_step}, we find a path in $T$ with vertices $ t_0, t_1, \ldots, t_s, t_{s+1}$ in order such that $t_0=t,t_1=t'$ and $t_it_{i+1}$ is $(2a(s-i)+3)$-good for $0 \leq i \leq s$. By \cref{lem:l_boundary} $|B_{t_i} \cap B_{t_{i+1}}| \geq n$  and so  $|B_{t_{i+1}} - B_{t_{i}}| \leq a$ for $0 \leq i \leq s$. Thus $|\cup_{i=0}^{s+1}B_{t_i}| \leq n + (s+2)a$ and $|B_{t_i}| \geq n$. By averaging there exists $v \in V(G)$ that belongs to at least $$\frac{(s+2)n} {n+(s+2)a} \geq \frac{n}{n/s + a} \geq \frac{n}{2n/s} = \frac{s}{2}$$ bags of $\mc{T}$, as claimed.		
\end{proof}

The most technical part of the proof is encapsulated in the following lemma.
	
	\begin{lemma}\label{lem:l_main} There exists $0 < \eps  < 1$ (independent of $n$) such that the following holds.
		Let $B \subseteq V(G)$ be such that $|B| \leq (1 + \eps)n$ and every component of $G \setminus B$ is non-global. Then there exists a $(k \times k)$-semigrid $X$ such that $X \subseteq G \setminus B$ and $k \geq \eps n$.
	\end{lemma}
	
As far as we know the lemma holds without the assumption that every component of $G \setminus B$ is non-global, but it is helpful in the proof and holds in the application.

Before proving \cref{lem:l_main}, let us derive \cref{thm:grid_linear_spread} from it and the already established results.

\begin{proof}[Proof of \cref{thm:grid_linear_spread} assuming \cref{lem:l_main}.] Let $\delta = \eps / 15$, where $\eps $ is as in \cref{lem:l_main}. If $\eps n \leq 15a$ then $\delta n /a  \leq 1$ the theorem trivially holds, and so we assume that $\eps n \geq 15a \geq 6a+9$.  By Observation~\ref{obs:separator} there exists $t \in V(T)$ such that $V_{tt'} - B_t$ is non-global for every $tt' \in E(T)$.  By \cref{lem:l_main} there exists a $(k \times k)$-semigrid $X$ such that $X \subseteq G \setminus B_t$ and $k \geq \eps n$. As $X$ is connected there exists $tt' \in E(T)$ such that $X \subseteq V_{tt'} - B_t$. It follows that $tt'$ is $k$-good. Let $$s = \left\lfloor \frac{\eps n - 3}{2a} \right\rfloor \geq  \frac{\eps n - 3}{2a} -1 \geq  \frac{\eps n }{3a },$$ where the last inequality holds by our lower bound on $n$.  By the choice of $s$ we have $2sa+3 \leq k$ and by Corollary~\ref{c:induct} the spread of $\mc{T}$ is at least $s/2 \geq  \frac{\eps n }{6a} > \frac{\delta {n}}{a}$, as needed.
 \end{proof}
 
It remains to prove \cref{lem:l_main}.

\begin{proof}[Proof of  \cref{lem:l_main}.] 
We assume that $n$ is sufficiently large as a function of $\eps$ to satisfy the inequalities throughout the proof. This can be done, as we can always reduce $\eps$ so that the lemma trivializes for $n$  outside of the assumed range.

With this additional assumption, we show that $\eps = 1/40$ satisfies the lemma. As every component of $G \setminus B$ is non-global, we assume without loss of generality that there exists no path in $G \setminus B$ with one end in $F$ and another in $C$.

Let $G'$ be any graph obtained from the $G$ by adding a diagonal edge in each square, that is, for each pair $i, j \in [n-1]$, we add exactly one of the edges $(i,j)(i+1,j+1)$ and $(i+1,j)(i,j+1)$. The following result is well-known (see e.g.~\cite{gale1979game} and the introduction of~\cite{hamkins2023infinite} for a sketch of the result below). 
\begin{theorem}[Variant of Hex theorem]
\label{thm:hex}
    For any set $B \subseteq V(G')$, either there is a path in $G'[B]$ from $L$ to $R$
    or there is a path in $G' \setminus B$ from $F$ to $C$.
\end{theorem}

\begin{cor}\label[cor]{hexvariant} There exists a graph $G'$, obtained from $G$ by adding diagonals to the squares of the grid, which contains a path $P$ in $G'[B]$ with one end in $L$ and another in $R$.\end{cor}

\begin{proof} In each square of $G$, we can add a diagonal edge while ensuring that $G'$ still does not have a path in $G' \setminus B$ from $F$ to $C$. In particular, if any corner of a cell is in $B$, we can add the diagonal between this corner and its opposite corner without creating any new paths contained entirely in $G' \setminus B$. If instead, all four corners of a cell are in $G \setminus B$, we can add either diagonal. There is already a path between any two corners contained entirely in $G \setminus B$, so if there were a path in $G' \setminus B$ with one end in $F$ and one end in $C$, there would have been a path in $G \setminus B$ with one end in $F$ and one end in $C$. With these diagonals added, by \cref{thm:hex}, there must be a path $P$ with one end in $L$ and another in $R$ such that $P$ only uses vertices of $B$. \end{proof}

By \cref{hexvariant}, there exists a sequence of pairwise distinct vertices $$(x_1,y_1),(x_2,y_2),\ldots,(x_m,y_m)$$ in $B$ such that $x_1=1, x_m=n$ and $|x_i - x_{i+1}|,|y_i -y_{i+1}| \leq 1$ for every $1 \leq i \leq m-1$. Clearly, $m \geq n-1$.
Let $k =\lceil \eps n \rceil$.
Let $B_0 = \{(x_i,y_i)\}_{i=1}^{m}$ and let $B_1= B - B_0$. Then $|B_1| \leq (1+\eps)n - m \leq k + 1.$

Let $y_{MIN} = \min_{1 \leq i \leq m} y_i$. If $y_{MIN} \geq 2k + 2$ then there exists $J \subseteq [n]$ with $\max J < y_{MIN}$ and $|J|=k$ such that none of the vertices in $B$ have second coordinate in $J$, and one can select $I \subseteq [n]$  with $|I|= k$ such that none of the vertices in $B_1$ have the first coordinate in $I$. Then $G \setminus B$ contains a straight $(k \times k)$-semigrid defined by the sets $I$ and $J$. 

Thus we assume $y_{MIN} \leq 2k + 1$. Similarly, defining $y_{MAX} = \max_{1 \leq i \leq m} y_i$ we assume
$y_{MAX} \geq n -2k-1$. Let $i',i'' \in [m]$   be such that $y_{i'}=y_{MIN}$ and $y_{i''}=y_{MAX}$. Then $|i' - i''| \geq y_{MAX} - y_{MIN} \geq n - 4k-2.$ As $m \leq n+k$, we may assume that $i' \leq 5k+2$ and $i'' \geq m - 5k-2$. It follows that $$y_1 \leq y_{MIN} + i'-1 \leq 7k+3 \qquad \mathrm{and} \qquad y_m \geq y_{MAX} - m + i'' \geq n - 7k -3.$$ 

We construct a sequence $1=i_1 < i_2 < \ldots < i_{\ell}=m$, recursively, as follows. Let $i_1=1$.  Once $i_1,\ldots,i_{j-1}$ are defined let $i_j \in [m-1]$ be chosen minimum such that $i_j > i_{j-1}$, $y_{i_{j}} > y_{j-1}$ and $x_{i_j} > x_{i_{j-1}}$ if such $i_j$  exists. If no such choice of $i_j$ is possible, set $\ell=j$, $i_{\ell}=m$ and terminate the construction. Let  $I_{*}=\{i_1,\ldots, i_{\ell}\}$. Our next goal is to show that $I_{*}$ is large.

For $i \in [m]$ let $s_i = x_i +y_i$.
Then for every $j \in [\ell-1]$ we have $x_{i_{j+1}} \leq x_{i_j}+1$ or $y_{i_{j+1}} \leq y_{i_j}+1$. Thus $s_{i_{j+1}} \leq s_{i_j} + (i_{j+1}-i_j)+1$, and so $s_{i_{\ell}} \leq s_1 + (i_{\ell}-i_1) + \ell -1$.
As $i_1=1$ and $i_{\ell}=m$, we have $s_m \leq s_1+ m+\ell-2$. But $s_1 \leq 7k+4$ and $s_m \geq 2n-7k-3$, implying \begin{equation} \label{e:ell} \ell \geq (2n-7k-3)+2 - 7k-4 - m \geq n-15k -5.\end{equation}

Let $B_* = \{(x_i,y_i)\}_{i \in I_* - \{m\}}$ and let $B_{**} = B - B_{*}$. Then $|B_{**}| \leq 16k + 6.$ 

For each $i \in [n]$ let $$U_i = \{(i,j) \: | \: i < j \} \cup \{(j , i) \: | \: j <i \}$$ and let $$D_i= \{(i,j) \: | \: i > j \} \cup \{(j , i) \: | \: j > i \}.$$ (See Figure~\ref{fig:tikz3} for an example.) Let $I_{**} = \{ i \in [n] \: | \: B_{**} \cap (U_i \cup D_i) \neq \emptyset\}.$ Then $|I_{**}| \leq 2|B_{**}| \leq 32k + 12$. Let $I_0 =  [n] - I_{**}$. 

We claim that for each $i \in I_0$ either $U_i$ or $D_i$ is disjoint from $B$. Indeed, both of these sets are disjoint from $B_{**}$. If $(x_j,y_j) \in B_* \cap U_i$ and $(x_{j'},y_{j'}) \in B_* \cap D_i$, then either $x_j = x_{j'}$ or $y_{j} = y_{j'}$, or $x_j < x_{j'}$ and $y_j >  y_{j'}$. None of these cases can occur by definition of $B_*$. Thus the claim holds.

Let $I_U = \{i \in I_0 \: | \: U_i \cap B = \emptyset\}$ and let $I_D = \{i \in I_0 \: | \: D_i \cap B = \emptyset\}$. Then $|I_U| + |I_D| \geq |I_0| \geq n - 32k - 12$. Thus we may assume by symmetry that $|I_U| \geq n/2 - 16k-6 \geq \s{\frac{1}{2}-16\eps}n-22 \geq \eps n,$  
where the last inequality holds by the choice of $\eps$. 
Thus $G \setminus B$ contains a crooked $(k \times k)$-semigrid corresponding to any subset of indices in $I_U$ of size $k$.

\begin{figure}[ht!]
	\centering

\begin{tikzpicture}[scale=0.9]
	% --- parameters ---
	\def\t{0.35}   % thickness / padding of the bands
	\def\e{0.18}   % half-size of the excluded neighborhood around (3,3); must satisfy \e < \t
	
	% --- styles ---
	\tikzset{
		vtx/.style={circle,draw,fill=white,inner sep=1.2pt},
		Ushape/.style={draw=green!50!black, line width=0.35pt, fill=green!30, fill opacity=0.25, draw opacity=0.7},
		Dshape/.style={draw=red!60!black,   line width=0.35pt, fill=red!30,   fill opacity=0.25, draw opacity=0.7}
	}
	
	% --- U_3 region (connected, but excluding a neighborhood of (3,3)) ---
	% horizontal arm: y=3, x=1..(3-e)
	\path[Ushape] (1-\t,3-\t) rectangle (3-\e,3+\t);
	% vertical arm: x=3, y=(3+e)..8
	\path[Ushape] (3-\t,3+\e) rectangle (3+\t,8+\t);
	% connector above (3,3), avoiding the excluded square
	\path[Ushape] (3-\e,3+\e) rectangle (3+\t,3+\t);
	
	% --- D_3 region (connected, but excluding a neighborhood of (3,3)) ---
	% vertical arm: x=3, y=1..(3-e)
	\path[Dshape] (3-\t,1-\t) rectangle (3+\t,3-\e);
	% horizontal arm: y=3, x=(3+e)..8
	\path[Dshape] (3+\e,3-\t) rectangle (8+\t,3+\t);
	% connector below (3,3), avoiding the excluded square
	\path[Dshape] (3+\e,3-\t) rectangle (3+\t,3-\e);
	
	% --- labels (black), moved as requested ---
	\node[black, anchor=south] at (3,8+\t+0.1) {$U_3$};
	\node[black, anchor=west]  at (8+\t+0.1,3) {$D_3$};
	
	% --- vertices ---
	\foreach \x in {1,...,8}{
		\foreach \y in {1,...,8}{
			\node (v\x\y) at (\x,\y) [vtx] {};
		}
	}
	
	% --- edges ---
	\foreach \x in {1,...,7}{
		\foreach \y in {1,...,8}{
			\pgfmathtruncatemacro{\xp}{\x+1}
			\draw (v\x\y) -- (v\xp\y);
		}
	}
	\foreach \x in {1,...,8}{
		\foreach \y in {1,...,7}{
			\pgfmathtruncatemacro{\yp}{\y+1}
			\draw (v\x\y) -- (v\x\yp);
		}
	}
\end{tikzpicture}

\caption{The sets $U_3$ and $D_3$ in an $(8 \times 8)$-grid.}
\label{fig:tikz3}
\end{figure}
\end{proof}

\section{Proof of Theorem~\ref{thm:average_spread_inf_is_1}}
We will use the following auxiliary result.
\begin{lemma}[Lemma 22 in \cite{wood2025treedecompositionssmallwidth}]
\label{lem:Wood_division}
For any integer $k \geq 2$, every rooted tree $T$ with $|V (T)|\geq k$ with root $r$ has a sequence $(T_1,\dots, T_m)$ of pairwise edge-disjoint rooted subtrees of $T$ such that:
\begin{itemize}
    \item $T=T_1\cup \dots \cup T_m$;
    \item $r\in V(T_1)$ and for $i \in \{2,\dots, m\}$, if $r_i$ is the root of $T_i$ then $V(T_i) \cap V (T_1 \cup \dots \cup T_{i-1})=\{r_i\}$;
    \item $|V(T_i)|\in \{k,\dots, 2k-2\}$ for each $i\in \{2,\dots,m\}$ and $|V(T_1)|\leq 2k-2$.
\end{itemize}
\end{lemma}

In fact, our proof of Theorem~\ref{thm:average_spread_inf_is_1} will follow an approach from Wood~\cite{wood2025treedecompositionssmallwidth}, for which we apply the lemma with a different parameter. Only the analysis of the spread in this construction is new.
\begin{proof}[Proof of Theorem~\ref{thm:average_spread_inf_is_1}]
Let $c'>1$. Choose an integer $t\geq 2$ such that
\[
t/(t-1)\leq c'.
\]
Let $c=3+4t$. We will prove that every graph $G$ has a
tree-decomposition $(T,(B_t)_{t\in V(T)})$ of width at most $c\tw(G)$ and average spread 
\[
\frac1{|V(G)|}\sum_{x\in V(T)}|B_x|\leq c'.
\]
Let $G$ be a graph of treewidth $w\geq 1$. Then there is a tree decomposition $(T,(B_x)_{x\in V(T)})$ of $G$ with width $w$, which we can root such that if $x$ is the parent of $y$ in $T$, then $|B_y\setminus B_x|=1$, that is, each bag contributes exactly one new vertex. We apply Lemma~\ref{lem:Wood_division} to $T$ with $k=t\cdot (w+1)+1$ to obtain $(T_1,\dots,T_m)$. (If the lemma does not apply, we are already done by putting all vertices in a single bag.)
Let \[
B_1'=\cup_{s\in V(T_1)}B_s
\]
and for $i\in \{2,\dots,m\}$, let
\[
B_i'=\cup_{s\in V(T_i)\setminus \{r_i\}}B_s.
\]
In \cite[Lemma 23]{wood2025treedecompositionssmallwidth} it is shown that such a \textit{quotient} yields a tree decomposition  $(F,(B_x')_{x\in V(F)})$ of $G$, for any tree $F$ with vertex-set $\{1,\dots, m\}$, rooted at vertex $1$, where for $i \in  \{2,\dots , m\}$, the parent
of $i$ is the smallest number $j\in \{1,\dots,i-1\}$ such that $r_i\in V(T_j)$.

Note that for any $i\in [m]$,
\[
|B_i'|\leq (w+1)+(|V(T_i)|-1)\leq w+2k \leq (1+4t+2)w=c\cdot \text{tw}(G),
\]
since moving from a parent to a child bag introduces \textit{at most} one new vertex.
Moreover, for all $i\in\{2,\dots,m\}$,  \[
|B_i'|\geq |V(T_i)|-1\geq k-1=t\cdot (w+1).
\]
since each bag introduces \textit{at least} one new vertex.

We say that vertex $v$ is \textit{new} in bag $B_i'$ if $i$ is the smallest $j\in [m]$ for which $v\in B_j'$. 
\begin{claim}
    For all $i\in [m]$, $B_i'$ has at least $(1-1/t)|B_i'|$ new vertices. 
\end{claim}
\begin{proof}
The new vertices in $B_i'$ are exactly those in $B_i'\setminus B_{r_i}$. Indeed,  $T_i$ only intersects $T_1,\dots,T_{i-1}$ in $r_i$ and since $(T,(B_s)_{s\in V(T)})$ is a tree decomposition,
\[
\left(\cup_{s\in V(T_i)\setminus \{r_i\}}B_s\right)\cap \left(\cup_{s\in V(T_1)\cup \dots \cup V(T_{i-1})}B_s\right)\subseteq B_{r_i}.
\]
Using that for each $s\in V(T_i)\setminus \{r_i\}$, there is at least one vertex in $B_s$ that was not in the bag of its parent, we find that
\[
|B_{r_i}\cap B_{i}'|/|B_{i}'|\leq (w+1)/(|V(T_i)|-1)\leq (w+1)/((w+1)t)= 1/t,
\]
and so there are at least $(1-1/t)|B_i'|$ new vertices in $B_i'$.
\end{proof}
We now show how the claim above implies that the average spread $\sum_{x\in V(F)}|B_x'|/|V(G)|$ of the tree decomposition is at most 
\[
t/(t-1)=1+1/(t-1)
\]
For $v\in V(G)$, let $i_v$ denote the (unique) index $i$ such that $v\in B_i'$ is new. Then
\[
|V(G)|=|\{(v,B_{i_v}'):v\in V(G)\}|=\sum_{i=1}^m |\{(v,B_i'):i=i_v,v\in V(G)\}|.
\]
By the claim, for each $i$, since the summand is exactly the number of new vertices,
\[
|\{(v,B_i'):i=i_v,v\in V(G)\}|\geq (1-1/t)|B_i'|.
\]
This proves
\[
|V(G)|\geq \sum_{i=1}^m |B_i'| (1-1/t).
\]
Rearranging shows that the average spread
\[
\sum_{i=1}^m |B_i'|/|V(G)|\leq 
1/(1-1/t) = 1/((t-1)/t)=t/(t-1).
\]
Hence we provided the desired tree decomposition.
\end{proof}

\section{Conclusion}
In this paper we increased our understanding of the trade-off between width and spread for tree decompositions. We discuss the algorithmic aspects and remaining open questions below.

\paragraph*{Algorithmic aspects.}
The proof of Theorem~\ref{thm:spread_constant_min_width_upper_bound2} can easily be made algorithmic.
Given a $c$-approximation algorithm for finding balanced separators, for any integer $b\geq 5$ a tree decomposition can be computed for any graph $G$ of width at most $c(3+8/b)(\tw(G)+1)$ for which each vertex $v\in V(G)$ is in at most $2b(d(v)+1)$ bags, at only polynomial overhead cost. Similarly, the proof of Theorem~\ref{thm:average_spread_inf_is_1} can also be made algorithmic, showing that for any $c'>1$, there is a constant $c>0$ and a polynomial time algorithm which computes for any graph $G$ a tree-decomposition of width at most $c\tw(G)$ and average spread at most $c'$.

Although notions related to spread have been motivated several times for potential algorithmic use, we are not currently aware of any algorithms that explicitly use tree decompositions with small spread and width. A recent application of a result on tree decompositions with additional properties, namely when the tree of the decomposition is relatively small, is in the construction of universal graphs~\cite{universalgraphs}.

\paragraph*{Tighter bounds on treewidth needed to bound the spread.}

A restriction of
\cref{thm:spread_constant_min_width_upper_bound2} to  graphs of maximum degree at most three, shows that for every $s \geq 8$ every  graph $G$ with $\Delta(G) \leq 3$  has a  tree decomposition of width at most $3\tw(G) + O(\tw(G)/s)$ and spread at most $s$. Conversely, \cref{lowermain} shows that for every $s$  there exist a graph $G$ of arbitrarily large treewidth  with $\Delta(G)=3$ such that every tree decomposition of $G$ of spread at most $s$ has width at least $3\tw(G) - 7$. We do not know which of these bounds on width is closer to the truth, i.e. the following two questions are open:
\begin{itemize}
\item
Does there exist $C > 0$ such that every graph $G$ with $\Delta(G) \leq 3$ admits a tree-decomposition of width at most  $3\tw(G) + C$ and spread at most $C$? 
\item Alternatively, does there exist  $\eps > 0$ such that for every $C,w > 0$ there exists a graph $G$ with $\Delta(G) \leq 3$ and $\tw(G) \geq w$ such that every tree decomposition of $G$ of spread at most $C$ has width at least $3\tw(G) + \eps \tw(G)/C$?
\end{itemize} 

\paragraph*{Declaration on AI use.} We used AI tools such as chatGPT for creating the TikZ figures and for proofreading. All proofs and writing were done entirely by the authors.

\bibliography{references}

\begin{thebibliography}{10}

\bibitem{Bienstock90}
Dan Bienstock.
\newblock On embedding graphs in trees.
\newblock {\em Journal of Combinatorial Theory, Series B}, 49(1):103--136, 1990.
\newblock \href {https://doi.org/10.1016/0095-8956(90)90066-9} {\path{doi:10.1016/0095-8956(90)90066-9}}.

\bibitem{bodlaenderlineartimealg}
Hans~L. Bodlaender.
\newblock A linear-time algorithm for finding tree-decompositions of small treewidth.
\newblock {\em SIAM Journal on Computing}, 25(6):1305--1317, 1996.
\newblock \href {https://doi.org/10.1137/S0097539793251219} {\path{doi:10.1137/S0097539793251219}}.

\bibitem{Bodlaender99}
Hans~L. Bodlaender.
\newblock A note on domino treewidth.
\newblock {\em Discrete Mathematics and Theoretical Computer Science}, 3(4):141--150, 1999.
\newblock \href {https://doi.org/10.46298/DMTCS.256} {\path{doi:10.46298/DMTCS.256}}.

\bibitem{Bodlaender24}
Hans~L. Bodlaender.
\newblock Approximation algorithms for treewidth, pathwidth, and treedepth --- {A} short survey.
\newblock In Daniel Kr{\'{a}}l and Martin Milanic, editors, {\em 50th International Workshop on Graph-Theoretic Concepts in Computer Science, {WG} 2024}, volume 14760 of {\em Lecture Notes in Computer Science}, pages 3--18. Springer, 2024.
\newblock \href {https://doi.org/10.1007/978-3-031-75409-8\_1} {\path{doi:10.1007/978-3-031-75409-8\_1}}.

\bibitem{BodlaenderE97}
Hans~L. Bodlaender and Joost Engelfriet.
\newblock Domino treewidth.
\newblock {\em Journal of Algorithms}, 24:94--123, 1997.
\newblock \href {https://doi.org/10.1006/jagm.1996.0854} {\path{doi:10.1006/jagm.1996.0854}}.

\bibitem{BodlaenderGHK95}
Hans~L. Bodlaender, John~R. Gilbert, H.~Hafsteinsson, and Ton Kloks.
\newblock Approximating treewidth, pathwidth, frontsize, and minimum elimination tree height.
\newblock {\em Journal of Algorithms}, 18:238--255, 1995.
\newblock \href {https://doi.org/10.1006/JAGM.1995.1009} {\path{doi:10.1006/JAGM.1995.1009}}.

\bibitem{BodlGroJacob}
Hans~L. Bodlaender, Carla Groenland, and Hugo Jacob.
\newblock On the parameterized complexity of computing tree-partitions.
\newblock {\em Discrete Mathematics \& Theoretical Computer Science}, 26:3, 2025.
\newblock \href {https://doi.org/10.46298/dmtcs.12540} {\path{doi:10.46298/dmtcs.12540}}.

\bibitem{BodlaenderGNS24}
Hans~L. Bodlaender, Carla Groenland, Jesper Nederlof, and C{\'{e}}line M.~F. Swennenhuis.
\newblock Parameterized problems complete for nondeterministic {FPT} time and logarithmic space.
\newblock {\em Information and Computation}, 300:105195, 2024.
\newblock \href {https://doi.org/10.1016/J.IC.2024.105195} {\path{doi:10.1016/J.IC.2024.105195}}.

\bibitem{COURCELLE199012}
Bruno Courcelle.
\newblock The monadic second-order logic of graphs. {I}. {R}ecognizable sets of finite graphs.
\newblock {\em Information and Computation}, 85(1):12--75, 1990.
\newblock \href {https://doi.org/10.1016/0890-5401(90)90043-H} {\path{doi:10.1016/0890-5401(90)90043-H}}.

\bibitem{DingO95}
Guoli Ding and Bogdan Oporowski.
\newblock Some results on tree decomposition of graphs.
\newblock {\em Journal of Graph Theory}, 20(4):481--499, 1995.
\newblock \href {https://doi.org/10.1002/JGT.3190200412} {\path{doi:10.1002/JGT.3190200412}}.

\bibitem{Distel2024}
Marc Distel and David~R. Wood.
\newblock Tree-partitions with bounded degree trees.
\newblock In David~R. Wood, Jan de~Gier, and Cheryl~E. Praeger, editors, {\em 2021--2022 MATRIX Annals}, volume~5 of {\em MATRIX Book Series}, pages 203--212. Springer, 2024.
\newblock \href {https://doi.org/10.1007/978-3-031-47417-0_11} {\path{doi:10.1007/978-3-031-47417-0_11}}.

\bibitem{DowneyM04}
Rodney~G. Downey and Catherine McCartin.
\newblock Online problems, pathwidth, and persistence.
\newblock In Rodney~G. Downey, Michael~R. Fellows, and Frank K. H.~A. Dehne, editors, {\em 1st International Workshop on Parameterized and Exact Computation, {IWPEC} 2004}, volume 3162 of {\em Lecture Notes in Computer Science}, pages 13--24. Springer, 2004.
\newblock \href {https://doi.org/10.1007/978-3-540-28639-4\_2} {\path{doi:10.1007/978-3-540-28639-4\_2}}.

\bibitem{DowneyM05}
Rodney~G. Downey and Catherine McCartin.
\newblock Bounded persistence pathwidth.
\newblock In Mike~D. Atkinson and Frank K. H.~A. Dehne, editors, {\em 11th Symposium on Computing: the Australasian Theory Symposium, {CATS} 2005}, volume~41 of {\em {CRPIT}}, pages 51--56. Australian Computer Society, 2005.

\bibitem{FominHT04}
Fedor~V. Fomin, Pinar Heggernes, and Jan~Arne Telle.
\newblock Graph searching, elimination trees, and a generalization of bandwidth.
\newblock {\em Algorithmica}, 41(2):73--87, 2005.
\newblock \href {https://doi.org/10.1007/S00453-004-1117-Y} {\path{doi:10.1007/S00453-004-1117-Y}}.

\bibitem{gale1979game}
David Gale.
\newblock The game of {H}ex and the {B}rouwer fixed-point theorem.
\newblock {\em The American Mathematical Monthly}, 86(10):818--827, 1979.
\newblock \href {https://doi.org/10.2307/2320146} {\path{doi:10.2307/2320146}}.

\bibitem{hamkins2023infinite}
Joel~David Hamkins and Davide Leonessi.
\newblock {Infinite Hex is a draw}.
\newblock {\em Integers: Electronic Journal of Combinatorial Number Theory}, 23(G3), 2023.
\newblock \href {https://doi.org/10.5281/zenodo.10075843} {\path{doi:10.5281/zenodo.10075843}}.

\bibitem{Jacob25}
Hugo Jacob, William Lochet, and Christophe Paul.
\newblock On a tree-based variant of bandwidth and forbidding simple topological minors.
\newblock {\em arXiv:2502.11674}, 2025.
\newblock \href {https://doi.org/10.48550/arXiv.2502.11674} {\path{doi:10.48550/arXiv.2502.11674}}.

\bibitem{universalgraphs}
Neel Kaul, Jaehoon Kim, Minseo Kim, and David Wood.
\newblock On universal graphs for trees and tree-like graphs.
\newblock {\em arXiv:2511.22358}, 2025.
\newblock URL: \url{https://arxiv.org/abs/2511.22358}.

\bibitem{RobertsonS13}
Neil Robertson and Paul~D. Seymour.
\newblock Graph minors. {XIII}. {T}he disjoint paths problem.
\newblock {\em Journal of Combinatorial Theory, Series B}, 63:65--110, 1995.
\newblock \href {https://doi.org/10.1006/JCTB.1995.1006} {\path{doi:10.1006/JCTB.1995.1006}}.

\bibitem{ROBERTSON198449}
Neil Robertson and P.D. Seymour.
\newblock {Graph minors. {III}. {P}lanar tree-width}.
\newblock {\em Journal of Combinatorial Theory, Series B}, 36(1):49--64, 1984.
\newblock \href {https://doi.org/10.1016/0095-8956(84)90013-3} {\path{doi:10.1016/0095-8956(84)90013-3}}.

\bibitem{Seese85}
Detlef Seese.
\newblock Tree-partite graphs and the complexity of algorithms.
\newblock In Lothar Budach, editor, {\em 5th International Conference on Fundamentals of Computation Theory, {FCT} 1985}, volume 199 of {\em Lecture Notes in Computer Science}, pages 412--421. Springer, 1985.
\newblock \href {https://doi.org/10.1007/BFB0028825} {\path{doi:10.1007/BFB0028825}}.

\bibitem{Wood09}
David~R. Wood.
\newblock On tree-partition-width.
\newblock {\em European Journal of Combinatorics}, 30(5):1245--1253, 2009.
\newblock \href {https://doi.org/10.1016/J.EJC.2008.11.010} {\path{doi:10.1016/J.EJC.2008.11.010}}.

\bibitem{wood2025treedecompositionssmallwidth}
David~R. Wood.
\newblock Tree decompositions with small width, spread, order and degree.
\newblock {\em arXiv:2509.01140}, 2025.
\newblock URL: \url{https://arxiv.org/abs/2509.01140}.

\end{thebibliography}

\appendix

\section{Omitted proof}
Recall that  Lemma~\ref{lem:good_spread_width_for_grids} states that for any $c>1$, any $(m\times n)$-grid has a path decomposition of width $cn$ and spread at most $\lceil 1/(c-1)\rceil+1$.
\begin{proof}[Proof of Lemma~\ref{lem:good_spread_width_for_grids}]
We divide $[n]$ into at most $b=\lceil1/(c-1)\rceil$ consecutive intervals of size at most $\lceil(c-1)n\rceil$, say $[n]=A_1\sqcup A_2\sqcup \dots \sqcup A_b$.
Let $A_{[i,j]}=A_i\cup A_{i+1}\cup \dots \cup A_{j}$. 
We claim the following sequence of bags forms the desired path decomposition.
The sequence starts with
\begin{align*}
(\{1\}\times A_{[1,b]})&\cup (\{2\}\times A_{1})\\
(\{1\}\times A_{[2,b]})&\cup (\{2\}\times A_{[1,2]})\\
(\{1\}\times A_{[3,b]})&\cup (\{2\}\times A_{[1,3]})\\
&\vdots \\
(\{1\}\times A_{b})&\cup (\{2\}\times A_{[1,b]})\\
(\{2\}\times A_{[1,b]})&\cup (\{3\}\times A_{1})\\
\end{align*}
and continues in this fashion.
To be precise, for $i\in [b]$ and $j\in [m-1]$, let
\[
B_{i,j}=(\{j\}\times A_{[i,b]})\cup (\{j+1\}\times A_{[1,i]})
\]
where $[1,1]=\{1\}$ and $[b,b]=\{b\}$, and let \[B_m = \{m\} \times [n] .\]
The sequence is 
\[
B_{1,1},B_{2,1},\dots,B_{b,1},B_{1,2},\dots,B_{b,2},\dots,B_{1,m-1},\dots,B_{b,m-1}, B_m.
\]
Note that each bag $B_{i,j}$ has size at most $n+\lceil(c-1)n\rceil=\lceil cn\rceil $, so the width is at most $cn$ as claimed.
Each vertex of the $m\times n$ grid is in $\{j\}\times A_i$ for exactly one value of $(i,j)$.
Each $\{j\}\times A_i$ appears in the sequence above at most $b+1$ times.
So the spread of each vertex is at most $b+1=\lceil1/(c-1)\rceil+1$.
Finally, we verify that the sequence of bags forms a path decomposition directly from the definition: the bags containing a vertex form a subpath and every vertex is covered.
All edges are of one of the following forms:
\[
\{(j,x),(j,x+1)\},~\{(j,x),(j+1,x)\}.
\]
Let $i$ be given such that $x\in A_i$.
Then $\{(j,x),(j,x+1)\} \in B_{i,j}$, and $\{(j,x),(j,x+1)\} \in B_{i,j}$,  for $j < m$, while $\{(m,x),(m,x+1)\} \in B_m$. So there is always a bag containing both types of edges above.
\end{proof}

\end{document}